\documentclass[11pt,a4paper,oneside]{amsart}

\usepackage{amssymb}
\usepackage{curves}
\input xy
\xyoption{all}

\newtheorem{definition}{Definition}
\newtheorem{thm}{Theorem}
\newtheorem{lem}{Lemma}
\newtheorem{prop}{Proposition}
\newtheorem{cor}{Corollary}
\newtheorem{rmk}{Remark}

\newcommand{\To}{\longrightarrow}
\newcommand{\Complex}{\mathbb C}
\newcommand{\Zeta}{\mathbb Z}
\newcommand{\Real}{\mathbb R}

\newcommand{\cplxi}{\mathbf i}
\newcommand{\D}[1]{\widetilde{D_#1}}
\newcommand{\E}[1]{\widetilde{E_#1}}
\newcommand{\cA}{\mathcal{A}}
\newcommand{\cF}{\mathcal{F}}
\newcommand{\cC}{\mathcal{C}}

\DeclareMathOperator{\End}{End}

\DeclareMathOperator{\tr}{tr}

\numberwithin{equation}{section}
\numberwithin{prop}{section}
\numberwithin{cor}{section}
\numberwithin{lem}{section}
\numberwithin{definition}{section}
\numberwithin{rmk}{section}

\begin{document}

\title[Certain Examples\dots]
{Certain Examples of Deformed Preprojective Algebras and Geometry
of Their *-Representations}

\author{Anton Mellit}

\email{mellit@mpim-bonn.mpg.de}

\address{Max Planck Institute for Mathematics, Vivatsgasse 7, D-53111 Bonn, Germany}


\begin{abstract}
  We consider algebras $e_i \Pi^\lambda(Q) e_i$ obtained from deformed preprojective algebra of
  affine type $\Pi^\lambda(Q)$ and an idempotent $e_i$ for certain concrete value of the vector $\lambda$ which corresponds to the traces of $-1\in SU(2, \Complex)$ in irreducible representations of finite subgroups of $SU(2, \Complex)$. We give a
  certain realization of these algebras which allows us to construct the
  $C^*$-enveloping algebras for them. Some well-known results, including description of four projections with sum 2 happen to be a particular case of this picture. 
\end{abstract}

\maketitle

\section{Introduction}
Let us take the *-algebra, generated by self-adjoint projections $a_1, a_2, a_3, a_4$ with relation $a_1+a_2+a_3+a_4 = 2$:
\[
A_{\D4}:=\Complex \langle a_1, a_2, a_3, a_4 | a_i=a_i^*, \sigma(a_i) \subset \{ 0, 1\}, a_1 + a_2 + a_3 + a_4 = 2 \rangle.
\]
Here and below $\sigma(a)$ denotes the spectrum of self-adjoint element $a$ and writing $\sigma(a)\subset \{x_1, x_2, \dots x_k\}$ for real numbers $x_1, x_2, \dots, x_k$ we mean that there is a relation 
\[
(a-x_1) (a-x_2) \dots (a-x_k) = 0.
\]
Classification of irreducible representations of $A_{\D4}$ is well known (see \cite{OS}). All irreducible representations, except some finite number of exceptional cases, form a two dimensional family of irreducible representations in dimension 2. The paper \cite{MyUmzh} was devoted to the classification of irreducible representations of the algebra
\[
A_{\E6}:=\Complex \langle a_1, a_2, a_3 | a_i=a_i^*, \sigma(a_i) \subset \{ 0, 1, 2 \}, a_1 + a_2 + a_3 = 3 \rangle.
\]
After that, the classification problem for irreducible representations of the *-algebra
\[
\begin{split}
A_{\E7}:=\Complex \langle a_1, a_2, a_3 | a_i=a_i^*, \sigma(a_1) \subset \{0, 1, 2, 3\}, \sigma(a_2) \subset \{0, 1, 2, 3\}, \sigma(a_3) \subset \{0, 2\}, \\ a_1 + a_2 + a_3 = 4 \rangle
\end{split}
\]
was solved in \cite{Ostrovskiy}. One can guess that the next step would be to write down formuli for irreducible representations of the *-algebra
\[
\begin{split}
A_{\E8}:=\Complex \langle a_1, a_2, a_3 | a_i=a_i^*, \sigma(a_1) \subset \{0, 1, 2, 3, 4, 5\}, \sigma(a_2) \subset \{0, 2, 4\}, \sigma(a_3) \subset \{0, 3\}, \\ a_1 + a_2 + a_3 = 6 \rangle.
\end{split}
\]
In this paper we treat four cases in a unified way. We obtain a description of the $C^*$-enveloping algebra for $A_{\D4}$, $A_{\E6}$, $A_{\E7}$, $A_{\E8}$. The $C^*$-enveloping algebra is represented as the algebra of continuous matrix-valued functions on a sphere which are equivariant with respect to a certain finite group which acts on the sphere and on matrices. Next few paragraphs contain the description of our approach to this problem.

Recall the definition of deformed preprojective algebra from \cite{CrB}.
Let $Q$ be a quiver with vertex set $I$. Write $\bar{Q}$ for the
double quiver of $Q$, i.e. the quiver obtained by adding a reverse arrow
$a^*: j\To i$ for each arrow $a: i\To j$, and write $\Complex \bar{Q}$
for its path algebra, which has basis the paths in $\bar{Q}$ including
a trivial path $e_i$ for each vertex $i\in I$. By weight we mean any
function from $I$ to complex numbers. Given weight $\lambda =
(\lambda_i)_{i\in I}$ the deformed preprojective algebra of weight
$\lambda$ is
\begin{equation}\label{PEq}
\Pi^\lambda(Q) = \Complex \bar{Q} / (\sum_{a\in Q}(a a^* - a^*a)
 -\sum_{i\in I} \lambda_i e_i),
\end{equation}

Suppose all $\lambda_i$, $i\in I$ are real. Then the
correspondence $e_i \To e_i$ for $i\in I$ and $a \To a^*$, $a^* \To a$
for $a\in Q$ can be uniquely extended to an antilinear involution of
$\Pi^\lambda(Q)$. We denote this involution by *.
 
Note that the
proposition 6.5 of \cite{King} can be reformulated in the following
way:
\begin{prop}
A representation $\phi$ of quiver $Q$ is a direct sum of
$\theta$-stable representations if and only if it can be extended to an
involutive representation of the deformed preprojective algebra
$\Pi^\theta(Q)$ with respect to the involution introduced above. If such
an extension exists it is unique up to an isomorphism of involutive
representations.
\end{prop}

Next we consider the decomposition of $\Pi^\lambda(Q)$ with respect to
orthogonal idempotents $e_i$, $i\in I$. In particular we obtain
algebras $e_i \Pi^\lambda(Q) e_i$. If $Q$ is a star quiver,
i.e. a quiver of type
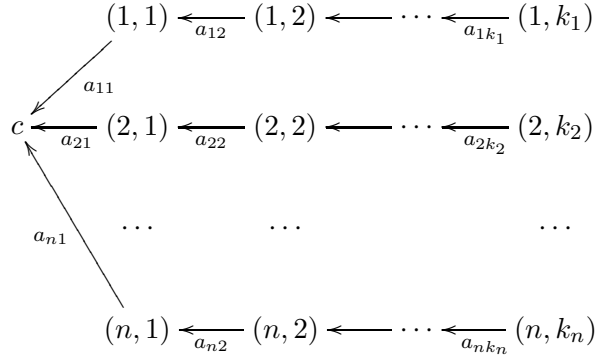
\[
        \xymatrix{
        & (1, 1) \ar@{->}[ld]^{a_{11}} & (1, 2) \ar@{->}[l]^{a_{12}} &
        \cdots \ar@{->}[l] & (1, k_1) \ar@{->}[l]^{a_{1 k_1}} \\
c       & (2, 1) \ar@{->}[l]^{a_{21}} & (2, 2) \ar@{->}[l]^{a_{22}} &
        \cdots \ar@{->}[l] & (2, k_2) \ar@{->}[l]^{a_{2 k_2}} \\
        & \cdots & \cdots & & \cdots\\
        & (n, 1) \ar@{->}[luu]^{a_{n1}} & (n, 2) \ar@{->}[l]^{a_{n2}}
        & \cdots \ar@{->}[l] & (n, k_n) \ar@{->}[l]^{a_{n k_n}} \\
        }
\]
then the algebra for the center vertex, choosing $x_i = a_{i1} a_{i1}^*$
can be described in terms of generators and relations as (see
\cite{CrB2}, \cite{MyWorkshop}):
\begin{equation} \label{ecPiec}
e_c \Pi^\lambda(Q) e_c = \Complex\langle x_1, \dots, x_n | P_i(x_i)=0\, (i=1,\dots,n),
    \sum_{i=1}^n x_i = \mu e\rangle,
\end{equation}
where $\mu = \lambda_c$ and 
\[
P_i(t) = (t-\alpha_{i 0}) (t-\alpha_{i 1}) \dots (t-\alpha_{i k_i}), \; \text{where} \;
\alpha_{i j} = -\sum_{l=1}^j \lambda_{i j}.
\]
With respect to the introduced involution these generators are
self-adjoint, so for an involutive representation $W$ the condition
$P_i(W_{x_i})=0$ is equivalent to
\[
\sigma(W_{x_i}) \subset \{\alpha_{i 0}, \alpha_{i 1}, \dots, \alpha_{i k_i} \},
\]
where $\sigma$ denotes the spectrum of a self-adjoint operator. Moreover, if numbers $\alpha_{i j}$ are pairwise distinct for each $i$ then $e_c \Pi^\lambda(Q) e_c$ is Morita equivalent to $\Pi^\lambda(Q)$, i. e. the functor from the category of representations of $\Pi^\lambda(Q)$ to the category of representations of $e_c \Pi^\lambda(Q) e_c$ given by $M \mapsto e_c M$ is an equivalence of categories. 

Note that replacing any $a\in Q$ by $a^*$ and $a^*$ by $-a$ does not change the relation (\ref{PEq}), so the deformed preprojective algebra does not depend on
the orientation of the quiver (although the involution depends) and statement (\ref{ecPiec}) remain valid for quivers of the following type:
\begin{equation}\label{quiver_type}
        \xymatrix{
        & (1, 1) \ar@{->}[ld]^{a_{11}} & (1, 2) \ar@{<-}[l]^{a_{12}} &
        \cdots \ar@{->}[l] & (1, k_1) \ar@{->}[l]^{a_{1 k_1}} \\
c       & (2, 1) \ar@{->}[l]^{a_{21}} & (2, 2) \ar@{<-}[l]^{a_{22}} &
        \cdots \ar@{->}[l] & (2, k_2) \ar@{->}[l]^{a_{2 k_2}} \\
        & \cdots & \cdots & & \cdots\\
        & (n, 1) \ar@{->}[luu]^{a_{n1}} & (n, 2) \ar@{<-}[l]^{a_{n2}}
        & \cdots \ar@{->}[l] & (n, k_n) \ar@{<-}[l]^{a_{n k_n}} \\
        }
\end{equation}
where the arrows in each ray are oriented in the alternating
way.

\begin{rmk} \label{rmk11}
Each *-representation in a Hilbert space of the deformed preprojective algebra $\Pi^\lambda(Q)$ with any orientation of the quiver $Q$ being restricted to the image of $e_c$ gives a *-representation of the algebra $e_c \Pi^\lambda(Q) e_c$. In the opposite direction the statement is wrong in general, but is true if the following inequalities hold for numbers $\alpha_{i j}$ for each $0\le i \le n$, $0\le j < k_i$, $j < l \le k_i$:
\begin{align*}
\alpha_{i j} < \alpha_{i l} &\;\text{if an arrow goes from $(i, j+1)$ to $(i, j)$ (we agree that $(i, 0)$ is $c$), or}\\
\alpha_{i j} > \alpha_{i l} &\;\text{if an arrow goes from $(i, j)$ to $(i, j+1)$.}
\end{align*}
This condition fixes certain ordering on the set of roots of each polynomial $P_i$, which depends on the orientation of the quiver.
\end{rmk}

\begin{rmk} \label{rmk12}
Replacing each $a$ by $a^*$ and $a^*$ by $a$ makes an isomorphism of *-algebras $\Pi^\lambda(Q)$ and $\Pi^{-\lambda}(Q')$ where $Q'$ is obtained from $Q$ by reversing all arrows.
\end{rmk}

For algebras $A_{\D4}$, $A_{\E6}$, $A_{\E7}$, $A_{\E8}$, the following quivers and weights give deformed preprojective algebras with equivalent category of representations:
\[
        \xymatrix{
        & 1 \ar@{<-}[d] & \\
1 \ar@{<-}[r] & -2 & 1 \ar@{<-}[l] \\
        & 1 \ar@{<-}[u] & \\
        }
       \xymatrix{
       & & 1 \ar@{<-}[d] & & \\
       & & -2 \ar@{->}[d] & & \\
1 \ar@{<-}[r] & -2 \ar@{->}[r] & 3 & -2 \ar@{->}[l] & 1 \ar@{<-}[l] \\
       }
\]
\[
       \xymatrix{
       & & & 2 \ar@{<-}[d] & & & \\
1 \ar@{<-}[r] & -2 \ar@{->}[r] & 3 \ar@{<-}[r] & -4 & 3 \ar@{<-}[l] &
       -2 \ar@{->}[l] & 1 \ar@{<-}[l] \\
       }
\]
\[
       \xymatrix{
       & & 3 \ar@{<-}[d] & & & & &\\
-2 \ar@{->}[r] & 4 \ar@{<-}[r] & -6 & 5 \ar@{<-}[l] &
       -4 \ar@{->}[l] & 3 \ar@{<-}[l] & -2 \ar@{->}[l] & 1 \ar@{<-}[l]
       \\
       }
\]
For example, the set $\{0, 1, 2, 3, 4, 5\}$ from the definition of $A_{\E8}$ must be ordered as $(0, 5, 1, 4, 2, 3)$, the set $\{0, 2, 4\}$  as $(0, 4, 2)$, the set $\{0, 3\}$ as $(0, 3)$. Taking differences one obtains sequences $(-5, 4, -3, 2, -1)$, $(-4, 2)$, $(-3)$. If we put them on a quiver of type (\ref{quiver_type}) the condition from the remark \ref{rmk11} will be satisfied. Applying the procedure from the remark \ref{rmk12} gives the quiver with numbers as on the picture.

Suppose $Q$ is an extended Dynkin quiver of type $\widetilde{A_n}$,
$\widetilde{D_n}$, $\widetilde{E_6}$, $\widetilde{E_7}$, or
$\widetilde{E_8}$.  Let $(\delta_i)_{i\in I}$ be the corresponding
minimal imaginary root and suppose $0\in I$ denotes the extending
vertex. For $\lambda : I \longrightarrow \mathbb C$ such that $\lambda
\cdot \delta = 0$ it is known that the deformed preprojective algebra
$\Pi^\lambda(Q)$ and all mentioned algebras $e_i \Pi^\lambda(Q) e_i$,
$i \in I$ have centers isomorphic to $e_0 \Pi^\lambda(Q) e_0$ and are
finitely generated modules over it (see \cite{CrB}, \cite{GE},
\cite{MyWorkshop}). The commutative ring $e_0 \Pi^\lambda(Q) e_0$ is
itself isomorphic to the coordinate ring of some fiber of the
semiuniversal deformation of the quotient singularity $\Complex^2 //
\Gamma$ where $\Gamma$ is the finite subgroup of $SU(2, \Complex)$
corresponding to $Q$ by the McKay correspondence. 

Recall that McKay correspondence assigns to each vertex $i$ of $Q$ an irreducible representation $V_i$ of $\Gamma$. An identity representation is assigned to the extending vertex $0$. $\delta_i$ gives the dimension of $V_i$. The number of edges between any $i$ and $j$~--- vertices of $Q$ equals to the number of times $V_i$ occurs in the decomposition of $V\otimes V_j$ into irreducibles. We denote by $V$ the tautological two-dimensional representation of $\Gamma$ as a subgroup of $SU(2, \Complex)$. We choose a hermitian structure on each $V_i$ which makes it an unitary representation.

Suppose $Q$ is bipartite (this includes all cases except
$\widetilde{A_n}$ with odd number of vertices), so that some vertices
are called odd and some are even. For the group $\Gamma$ it means that $\Gamma$ contains negative identity of $SU(2, \Complex)$. Let $0$ be even. Suppose all arrows are directed from
odd vertices to even ones, and $\lambda_i = - \delta_i$ for
odd vertices and $\lambda_i = \delta_i$ for even ones. Clearly $\lambda_i$ equals to the trace of $-I \in \Gamma$ in $V_i$. Let $i\in
I$. Consider $\End_{\mathbb C}(V_i)$ equipped with $\Gamma$-action by
conjugation. Let $S_\Real$ be the real unit sphere in $\Real^3$ and
$S_\Complex$ be the affine variety consisting of points $(x,
y, z)\in \mathbb C^3$, which satisfy $x^2+y^2+z^2=1$. To define
$\Gamma$-action on $S_\Real$ and $S_\Complex$ we use the
well-known homomorphism from $\Gamma \subset SU(2, \Complex)$ to
$SO(3, \Real) \subset SO(3, \Complex)$. The image of $\Gamma$ in $SO(3, \Real)$ is isomorphic to $\Gamma/\{I, -I\}$ which we denote by $\Gamma'$. Note that $\Gamma'$ acts on $\End_{\mathbb C}(V_i)$. We prove the following:
\begin{thm} \label{thm1}
  The algebra $e_i \Pi^\lambda(Q) e_i$ is isomorphic to the
  algebra of polynomial $\Gamma'$-equivariant  maps for $S_{\mathbb C}$
  to $\End_{\mathbb C}(V_i)$. The involution on the latter algebra given by the formula $f^*(x) = f(\overline{x})^*$, $x\in S_\Complex$, coincides with the one induced from the former algebra.
\end{thm}

If, furthermore, $V_i$ is not exceptional in the sense of the definition $\ref{def51}$ then it is possible to construct a $C^*$-enveloping algebra
\begin{thm}
If $V_i$ is not exceptional in the sense of the definition $\ref{def51}$ then the $C^*$-enveloping algebra of $e_i \Pi^\lambda(Q) e_i$ exists and is isomorphic to the
$C^*$-algebra of continuous $\Gamma'$-equivariant maps for $S_\Real$ to $\End_{\mathbb C}(V_i)$.
\end{thm}
Exceptional cases are listed in the section \ref{sec6}.

Suppose $V_i$ is not exceptional. We can see at this point that the $C^*$-envelope of $e_i
\Pi^\lambda(Q) e_i$ defines a bundle of algebras on $S_\Real /
\Gamma'$, which is homeomorphic to $S_\Real$. It occurs that this
bundle can be trivialized in the following sence.
Let $x\in S_\Real /
\Gamma'$ be any orbit in 
$S_\Real$. Choose some representative $x_0$ for $x$ in the fundamental region and consider
the stabilizer subgroup of $\Gamma'$ for $x_0$. We denote the set of all elements
in $\End_{\mathbb C}(V_i)$ which commute with the stabilizer by
$M_{x_0}$. Then, the *-algebra $M_{x_0}$ has the form
\[
M_{x_0} \cong \End(\Complex^{d_1}) \times \End(\Complex^{d_2}) \times \dots \times \End(\Complex^{d_k}) \subset \End(\Complex^d),
\]
with $d_1 \ge d_2 \ge \dots \ge d_k$ and $d$ is the dimension of $V_i$. Denote by $N_x$ the subalgebra of $\End(\Complex^d)$ given by the righthand side of the expression above.
 Clearly if the stabilizer lies in the center of $\Gamma'$ then
$M_x$ equals to the whole $\End_{\mathbb C}(V_i)$. It is easy to see
that this holds for all but finitely many points of $S_\Real /
\Gamma'$. The number of orbits for which the stabilizer does not lie in
the center is $2$ for $A_n$ graphs and $3$ for $D_n$ and $E_n$ graphs.
\begin{thm}
The $C^*$-envelope of $e_i \Pi^\lambda(Q) e_i$ is isomorphic to the
$C^*$-algebra of continuous functions $f$ from $S_\Real / \Gamma$ (which is homeomorphic to sphere) to
$\End(\Complex^d)$ such that for any $x\in S_\Real / \Gamma$
$f(x) \in N_x$.
\end{thm}

One can apply theorems 1-3 to the case of algebras $A_{\D4}$, $A_{\E6}$, $A_{\E7}$, $A_{\E8}$ (the corresponding representation of $\Gamma$ is not exceptional). Then $\Gamma'$ is, correspondingly, dihedral, tetrahedral, octahedral and icosahedral group, $\Gamma$ is its preimage in $SU(2, \Complex)$, $i=c$ and $V_i$ is the unitary irreducible representation of $\Gamma$ of maximal dimension. The dimension is, correspondingly, $2$, $3$, $4$ and $6$. For example, for $A_{\E8}$ the corresponding group $\Gamma'$ is the group of symmetries of the icosahedron. The stabilizers are non-trivial for centers of faces, centers of edges and vertices. The corresponding stabilizers are cyclic groups of order $3$, $2$ and $5$. The representation of $\Gamma$ $V_c$ is given by the maps
\[
\Gamma \hookrightarrow SU(2, \Complex) \longrightarrow SU(6, \Complex),
\]
where the last arrow is given by the symmetric $5$-th power. The theorem 3 together with studying stabilizers of vertices gives
\begin{cor}
The $C^*$ -envelope of $A_{\E8}$ is isomorphic to the $C^*$-algebra of continuous maps $f$ from the 2-sphere with three marked points $a$, $b$, $c$ to the $C^*$-algebra of $6\times 6$ matrices such that the matrix $f(a)$ has all zero entries except three $2\times 2$ blocks on the diagonal, the matrix $f(b)$ has all zero entries except two $3\times 3$ blocks on the diagonal and the matrix $f(c)$ has all zero entries except four diagonal entries and one $2\times 2$ block. It follows that the set of irreducible representations contains exactly $4$ one-dimensional representations corresponding to point $c$, $4$ two-dimensional representations~--- $3$ for $a$ and $1$ for $c$, $2$ three dimensional representations for $b$ and a family of six-dimensional representations parametrized by points of the sphere excluding $a$, $b$, $c$.
\end{cor}

Theorem 1 is proved in the section \ref{sec4} after reviewing the connection between deformed preprojective algebras $\Pi$ and certain skew group rings $S$ in the section \ref{sec3}. The $C^*$ enveloping algebra is studied in the section \ref{sec5}. The section \ref{sec6} is devoted to study of exceptional representations and the section \ref{sec7} contains the proof of the theorem 3.

\section{Notation}
By superalgebra we mean a $\Zeta/ 2 \Zeta$-graded associative algebra over $\Complex$. If $A$ is a superalgebra we denote by $A_{ev}$ the space in grade $0$ and by $A_{odd}$ the space in grade $1$ so that $A = A_{ev} \oplus A_{odd}$. Every algebra $A$ then can be considered as a superalgebra with $A_{odd}=0$.
\begin{definition} The involution on a superalgebra $A$ is a map $\circ: A \To A$ such that for any $a, b \in A$, $t\in \Complex$
\begin{enumerate}
\item $(a+b)^\circ = a^\circ + b^\circ$,
\item $(ab)^\circ = b^\circ a^\circ$,
\item $t^\circ = \bar{t}$,
\item $a^{\circ \circ} = a$ if $a\in A_{ev}$,
\item $a^{\circ \circ} = -a$ is $a\in A_{odd}$.
\end{enumerate}
The classical involution on an algebra $A$ is a map $\circ$ which satisfies properties (1)-(3) and $a^{\circ \circ} = a$ for all $a\in A$.
\end{definition} 

Let $G$ be an affine Dynkin graph of type ADE such that $G$ is not a cycle with odd number of vertices. Denote its vertex set by $I$ and let $0\in I$ be the extending vertex. Let $\delta=(\delta_i)$ be the minimal imaginary root for the root system corresponding to $G$. It is possible to split $I$ into two parts called even and odd in such a way that $0$ is even and every edge of $G$ connects one even and one odd vertex. We put $\sigma_i = -1$ if $i$ is odd and $\sigma_i = 1$ if $i$ is even.

Let $Q$ be a quiver obtained from $G$ by directing every edge towards an even vertex. Let $\lambda = (\lambda_i)$ be given by $\lambda_i = \sigma_i \delta_i$. 

By the McKay correspondence there exists a finite group $\Gamma \subset SU(2, \Complex)$ acting on $V=\Complex^2$ and a bijection $i \longleftrightarrow V_i$ between $I$ and the set of all nonisomorphic irreducible representations for $\Gamma$ such that
\begin{enumerate}
\item $V_0$ is the trivial representation.
\item $V_i \otimes V$ is isomorphic to the direct sum of $V_j$ where $j$ ranges over vertices of $I$ connected with $i$ by an edge.
\item $\dim V_i = \delta_i$.
\end{enumerate}
Note that for the cases under consideration $\Gamma$ contains an element $\tau = -I_V$ which belongs to the center of $\Gamma$ and $\tr_{V_i} \tau = \lambda_i$ for all $i\in I$.

Consider a skew group algebra $S = \Complex \langle V \rangle \ast \Gamma$ where $\Complex \langle V \rangle$ denotes the tensor algebra of $V^*$. Denote by $x$, $y$ the elements of the standard basis of $V^* = \Complex^{2*}$, by $\varepsilon_x$, $\varepsilon_y$ the elements of the dual basis of $V$ and by $w$ the element of $(V\otimes V)^* \cong V^* \otimes V^* \subset S$ given by $x\otimes y - y \otimes x$ so that $w(\varepsilon_x, \varepsilon_y) = 1$, $w(\varepsilon_y, \varepsilon_x) = -1$ and $w(\varepsilon_x, \varepsilon_x) = w(\varepsilon_y, \varepsilon_y) = 0$. Then $S$ is generated by $x$, $y$ and elements of $\Gamma$. $\Zeta$-grading of $\Complex \langle V \rangle$ induces a $\Zeta$-grading of $S$, which in its turn induces a $\Zeta/2 \Zeta$-grading. So $S$ is a superalgebra.

There exists a unique involution $\circ$ on $S$ such that
\begin{enumerate}
\item $g^\circ = g^{-1}$ for $g\in \Gamma$,
\item $(\cdot, v)^\circ = \cplxi w(\cdot, v)$ for $v\in V$,
\end{enumerate}
here $\cplxi = \sqrt{-1}$. One can calculate
\[
x^\circ = (\cdot, \varepsilon_x)^\circ = \cplxi w(\cdot, \varepsilon_x) = 
-\cplxi y,
\]
\[
y^\circ = (\cdot, \varepsilon_y)^\circ = \cplxi w(\cdot, \varepsilon_y) = \cplxi x.
\]
It follows that
\[
w = \cplxi x x^\circ + \cplxi y y^\circ\; \text{so} \; w^\circ = w.
\]

Take a factor algebra $S^\lambda = S/(xy - yx - \tau)$. Clearly, it is again a superalgebra and the involution $\circ$ of $S$ induces an involution on $S^\lambda$ which we will denote again by $\circ$.

Next, consider the path algebra of the double of $Q$ denoted by $\Pi$. It is generated by idempotents $e_i$ for each vertex $i\in I$, arrows $a\in Q$ and opposite arrows $a^*$ for $a\in Q$. This algebra is a superalgebra in an obvious way and we define an involution by
\begin{enumerate}
\item $e_i^\circ = e_i$ for $i\in I$,
\item $a^\circ = -\cplxi a^*$ for $a\in Q$,
\item $a^{*\circ} = \cplxi a$.
\end{enumerate}

The factor algebra $\Pi^\lambda$ is defined by
\[
\Pi^\lambda = \Pi / (\sum_{a\in Q}(a a^* - a^*a)
 -\sum_{i\in I} \lambda_i e_i),
\]
and is again a superalgebra with an induced involution.

There is also a classical involution $*$ on $\Pi$ which induces a classical involution on $\Pi^\lambda$. The action of $*$ is given by
\[
e_i^* = e_i, \; (a)^* = a^*, \; (a^*)^* = a.
\]

\section{Connection between $\Pi$ and $S$} \label{sec3}

Following \cite{CrB} choose an idempotent $f_i$ for each $i\in I$ such that $V_i \cong \Complex \Gamma f_i$. We additionally require $f_i$ to be self-adjoint and the Hermitian structure on $V_i$ to be induced from that of $\Complex \Gamma$. Put $f = \sum_{i\in I} f_i$. Then
\begin{prop}[see \cite{CrB, MyWorkshop}] \label{iso1}
There exists a graded isomorphism $\phi_1: f S f \cong \Pi$ (with respect to $\Zeta$-grading) such that
\begin{enumerate}
\item $\phi_1(f_i) = e_i$ for $i\in I$,
\item $\phi_1(\sum_{i\in I} \delta_i f_i(x y - y x) f_i) = \sum_{a\in Q} (a a^* - a^* a)$.
\end{enumerate}
This induces a graded isomorphism $f S^\lambda f \cong \Pi^\lambda$ (with respect to $\Zeta / 2 \Zeta$-grading).
\end{prop}
 
We prove a stronger result:
\begin{prop} \label{Prop32}
There exists a graded isomorphism $\phi: f S f \cong \Pi$ (with respect to $\Zeta$-grading) such that
\begin{enumerate}
\item $\phi$ satisfies conditions of the proposition \ref{iso1} and
\item $\phi(a^\circ) = \phi(a)^\circ$ for any $a\in f S f$.
\end{enumerate}
This induces a graded isomorphism $f S^\lambda f \cong \Pi^\lambda$ (with respect to $\Zeta / 2 \Zeta$-grading) which respects involutions.
\end{prop}
\begin{proof}
Suppose a collection of positive real numbers $(c_a)_{a\in Q}$ is given. Than there exists a unique graded automorphism $\phi_2: \Pi \To \Pi$ such that
\begin{enumerate}
\item $\phi_2(e_i) = e_i$ for $i\in I$,
\item $\phi_2(a) = c_a a$ for $a\in Q$,
\item $\phi_2(a^*) = c_a^{-1} a^*$ for $a\in Q$.
\end{enumerate}
We put $\phi = \phi_2 \circ \phi_1$. Clearly the first condition is satisfied. We need to prove that it is posible to choose numbers $c_a$ in such a way that the second condition would be satisfied. We denote the involution on $\Pi$ induced from $f S f$ by $\phi_2$ as $\star$.

Let $a\in Q$. Then $a^\star = t_a a^*$, $t_a\in \Complex$ because $\star$
preserves $\mathbb{Z}$-grading and there is at most one arrow
between each two vertices. Then $a^{\star\star}=-a$ implies
$a^{*\star} = -\bar{t_a^{-1}} a$ and since $w$ is self-adjoint
\[
\sum_{a\in Q} (a a^* - a^* a) = (\sum_{a\in Q} (a a^* - a^* a))^\star
= \sum_{a\in Q} -\frac{t_a}{\bar{t_a}} (a a^* - a^* a) \;
\text{implies} 
\]
\[
-\frac{t_a}{\bar{t_a}} = 1, \; a\in Q.
\]
Therefore $t_a = r_a \cplxi$ for some $r_a \in \Real$. Then we can
express
\begin{equation}\label{eq1}
\sum_{a\in Q} (a a^* - a^* a) = \sum_{a\in Q} ( -\frac{\cplxi}{r_a} a
a^\star - \cplxi r_a a^* a^{*\star}) = -\cplxi \sum_{a\in \bar{Q}} q_a
a a^\star
\end{equation}
for some real numbers $q_a$, $a\in \bar{Q}$. On the other hand 
\[
\sum_{i\in I} \delta_i f_i (x y - y x) f_i = \cplxi \sum_{i\in I}
\delta_i f_i (x x^\circ + y y^\circ) f_i.
\]
Consider the following element of the group algebra 
\[
J = \frac1{|\Gamma|} \sum_{g\in \Gamma} g f g^{-1}.
\]
Clearly it belongs to the center of the group algebra and its trace on
each irreducible representation equals to $1$. So $J^{-1}$ is central
and positive, hence there exists central self-adjoint $J'$ such that
$J' J J' = 1_{\Complex\Gamma}$. It follows that $1_{\Complex\Gamma}$
can be represented as 
\[
1_{\Complex\Gamma} = \sum_{k=1}^K \alpha_k f \alpha_k^\circ, \;
\alpha_k\in \Complex\Gamma.
\]
Using this we can represent $f_i x x^\circ f_i$ as
\[
f_i x x^\circ f_i = \sum_{k=1}^K (f_i x \alpha_k f \alpha_k^\circ x^\circ f_i) = \sum_{k=1}^K \sum_{j\in I} (f_i x \alpha_k f_j) (f_i x \alpha_k f_j)^\circ,
\]
which is mapped by $\phi_1$ to a linear combination with positive coefficients of elements of the form $a a^\star$ for $a\in\bar{Q}$. Applying the same
arguments for $f_i y y^\circ f_i$ we obtain that the sum
\[
\sum_{i\in I} \delta_i f_i (x x^\circ + y y^\circ) f_i 
\]
is mapped by $\phi_1$ to a linear combination with positive coefficients of elements of the form $a a^\star$ for $a\in\bar{Q}$. Hence the numbers $q_a$ in \ref{eq1} are negative and it follows that the numbers $r_a$ are negative too. 
We put $c_a = \sqrt{-r_a}$ for $a\in Q$ and obtain
\[
\phi_2(a^\star) = \cplxi r_a c_a^{-1} a^* = - \cplxi c_a a^* = \phi_2(a)^\circ \; \text{and}
\]
\[
\phi_2(a^{*\star}) = -\cplxi r_a^{-1} c_a a = \cplxi c_a^{-1} a = \phi_2(a^*)^\circ,
\]
so the composition $\phi_2 \circ \phi_1$ satisfies the second condition.
\end{proof}

\section{The even part} \label{sec4}
We are going to consider algebras $e_i \Pi^\lambda e_i$. Since the graph is bipartite we can replace $\Pi^\lambda$ by its even part, i.e. the following is clear
\begin{prop} \label{Prop41}
The superalgebras $e_i \Pi e_i$, $e_i \Pi^\lambda e_i$ have zero odd part. So
$e_i \Pi e_i = e_i \Pi_{ev} e_i$ and $e_i \Pi^\lambda e_i = e_i \Pi^\lambda_{ev} e_i$.
\end{prop}

Recall that the algebra $\Pi$ has a classical involution $*$.
\begin{prop} \label{Prop42}
The restrictions of $*$ and $\circ$ to $\Pi_{ev}$ coincide.
\end{prop}
\begin{proof}
It is enough to check the statement for elements of degree $2$. Since each arrow of the quiver goes from an odd vertex to an even one we have that products
$a b$ and $a^* b^*$ are zero for $a, b \in Q$. For $a^* b$ and $a b^*$ 
\[
(a^* b)^\circ = b^* a = (a^* b)^* \; \text{and} \; 
(a b^*)^\circ = b a^* = (a b^*)^*.
\]
\end{proof}

Consider the space $W$ of traceless operators on $V$. It is a complex $3$-dimensional vector space with symmetric bilinear form $g$ given by
\[
g(a, b) = 2 \tr(a b).
\]
Clearly the subspace $W_\Real$ of traceless hermitian operators on $V$ is a real $3$-dimensional vector space with scalar product induced by $g$. The group $SU(V)$ acts on $W$ and $W_\Real$ fixing the form $g$. This produces a well-known homomorphism
\[
SU(2, \Complex) \cong SU(V) \To SO(W_\Real) \cong SO(3, \Real)
\]
whose kernel is $\{-1, 1\}$. The subvariety $W_1\subset W$ of operators $a$ such that $g(a, a) = 1$ is again equipped with $\Gamma$-action.

Next, we identify $V\otimes V$ with $V\otimes V^*=
\End_\Complex(V)$ by the map $\psi$ defined as
\[
\psi(\alpha \otimes \beta) = \alpha w(\beta, \cdot), \; \alpha, \beta
  \in V.
\]
Clearly $\psi$ is equivariant and $\psi^{-1}$ restricts to an embedding $W_1 \To V\otimes V$ which induces an epimorphism of coordinate rings $\Complex[V\otimes V] \To \Complex[W_1]$ which in its turn induces an epimorphism
\[
\Psi_1: \Complex \langle V \rangle_{ev} \To \Complex[W_1]
\]
which is equivariant. Denote its kernel by $K$. Let $a = x^2$, $b = y^2$, $c = x y + y x$ and $d = x y - y x$ be generators of $\Complex \langle V \rangle_{ev}$.

\begin{prop} \label{kgen}
The ideal $K$ has the following set of generators:
\begin{equation} \label{gen1}
\{a b - b a, a c - c a, b c - c b, d, 4 a b - c^2 + 1\}
\end{equation}
\end{prop}
\begin{proof}
The kernel of the epimorphism $\Complex \langle V \rangle_{ev} \To \Complex[V \otimes V]$ is generated by pairwise commutators of generators $a$, $b$, $c$, $d$. Introduce a coordinates $m_{ij}$ ($i, j = 1, 2$) in $W$ such that for $m\in W$
\[
m = \begin{pmatrix}m_{11}(m) & m_{12}(m) \\ m_{21}(m) & m_{22}(m) \end{pmatrix}.
\]
Then for $t = \alpha \otimes \beta \in V \otimes V$
\[
\psi(t) = \alpha w(\beta, \cdot) = 
\begin{pmatrix} x(\alpha) \\ y(\alpha) \end{pmatrix}
\begin{pmatrix} -y(\beta) & x(\beta) \end{pmatrix} =
\begin{pmatrix} -\frac{c+d}2(t) & a(t) \\ -b(t) & \frac{c-d}2(t) \end{pmatrix}.
\]
Thus, taking the inverse map we obtain
\[
\Psi_1(a) = m_{12}, \; \Psi_1(b) = -m_{21}, \; \Psi_1(c) = m_{22} - m_{11},
\; \Psi_1(d) = - m_{11} - m_{22}.
\]
The equations of $W_1$ are
\[
m_{11} + m_{22} = 0, \; 2 m_{11}^2 + 4 m_{12} m_{21} + 2 m_{22}^2 = 1,
\]
so $K$ is generated by pairwise commutators of $a$, $b$, $c$, $d$ and two more elements:
\[
d \; \text{and} \; \frac{(c+d)^2}2 - 4 a b + \frac{(c-d)^2}2 - 1 = 
c^2 + d^2 - 4 a b - 1,
\]
and it can be easily seen that the set \ref{gen1} generates the same ideal.
\end{proof}

Denote by $I$ the ideal in $S$ generated by $x y - y x - \tau$ and by $I_{ev}$ the intersection $I \cap S_{ev}$. The following is true:
\begin{prop} \label{igen}
The ideal $I_{ev}$ is generated by
\[
a b - b a, a c - c a, b c - c b, d - \tau, 4 a b - c^2 + 1.
\]
\end{prop}
\begin{proof}
The ideal $I_{ev}$ can be generated by $x y - y x - \tau$ and $V^* (x y - y x - \tau) V^*$, so by five elements $g_1$, $g_2$, $g_3$, $g_4$, $g_5$, where
\begin{equation}
\begin{split}
x y - y x - \tau = d - \tau = g_1 \\
2 x (x y - y x - \tau) x = a (c-d) - (c+d) a + 2 a \tau
\equiv_{\mod g_1} a c - c a = g_2 \\
2 y (x y - y x - \tau) y = (c-d) b - b (c+d) + 2 b \tau
\equiv_{\mod g_1} c b - b c = g_3 \\
4 x (x y - y x - \tau) y = 4 a b - (c+d) (c+d) + 2 (c+d) \tau
\equiv_{\mod g_1} 4 a b - c^2 + 1 = g_4\\
4 y (x y - y x - \tau) x = (c-d) (c-d) - 4 b a + 2 (c-d) \tau
\equiv_{\mod g_1} c^2 - 4 b a - 1 = \\
\equiv_{\mod g_4} 4(a b - b a) = 4 g_5
\end{split}
\end{equation}
\end{proof}

Take an automorphism $\Psi_2$ of $S_{ev}$ which satisfies
\[
\Psi_2(g) = g\; (g\in\Gamma), \; \Psi_2(a) = a, \; \Psi_2(b) = b, \;
\Psi_2(c) = c, \; \Psi_2(d) = d+\tau.
\]
Since both $d$ and $\tau$ commute with elements of $\Gamma$ such an automorphism exists and propositions \ref{kgen} and \ref{igen} imply
\begin{cor} \label{Cor41}
The epimorphism $\Psi_0: S_{ev} \To \Complex[W_1] \ast \Gamma$ defined as $\Psi_0 = (\Psi_1 \otimes Id_{\Complex \Gamma}) \circ \Psi_2$ has kernel $I_{ev}$. Thus $\Psi_0$ induces an isomorphism $\Psi: S^\lambda_{ev} \cong \Complex[W_1] \ast \Gamma$.
\end{cor}

Define a classical involution $*$ on $\Complex[W_1]$ by formula
\[
(f^*)(m) = \overline{f(m^*)}.
\]
It induces a classical involution $*$ on $\Complex[W_1] \ast \Gamma$ since the operation of taking the hermitian adjoint commutes with the group action.
It occurs that 
\begin{cor} \label{Cor42}
The homomorphism $\Psi$ respects involutions.
\end{cor}
\begin{proof}
It is enough to prove the statement for homomorphisms $\Psi_1$ and $\Psi_2$.
Compute
\[
a^\circ = -b, \; b^\circ = -a, \; c^\circ = c, \; d^\circ = d, \; \tau^\circ = \tau,
\]
and
\[
m_{ij}^* = m_{ji} \; (i, j = 1, 2).
\]
Using defining formuli for $\Psi_1$ and $\Psi_2$ we prove the statement.
\end{proof}

\begin{cor}
There is an isomorphism $e_i \Pi^\lambda e_i \cong f_i \Complex[W_1] \ast \Gamma f_i$ which respects involutions.
\end{cor}
\begin{proof}
An isomorphism can be constructed using proposition \ref{Prop32} and corollary \ref{Cor41} taking into account proposition \ref{Prop41}. The fact that this isomorphism respects involution follows from propositions \ref{Prop42}, \ref{Prop32} and corollary \ref{Cor42}. All involved algebras can be displayed on the following diagram:
\[
\begin{array}{ccccccc}
e_i \Pi^\lambda e_i & \subset & \Pi^\lambda_{ev} & \twoheadleftarrow & \Pi_{ev} & \hookrightarrow & \Pi \\
\|\wr & & \|\wr & & \|\wr & & \|\wr \\
f_i S^\lambda f_i & \subset & f S^\lambda_{ev} f & \twoheadleftarrow & f S_{ev} f & \hookrightarrow & f S f \\
\|\wr & & \cap & & \cap & & \cap \\
f_i S^\lambda_{ev} f_i & \subset & S^\lambda_{ev} & \twoheadleftarrow & S_{ev} & \hookrightarrow & S \\
\|\wr & & \|\wr & & \|\wr & & \|\wr \\
f_i \Complex[W_1] \ast \Gamma f_i & \subset & \Complex[W_1] \ast \Gamma & \twoheadleftarrow & \Complex \langle V \otimes V \rangle \ast \Gamma & \hookrightarrow & \Complex \langle V \rangle \ast \Gamma 
\end{array}
\]
On the diagram relation $A \subset B$ means that the algebra $A$ is contained in $B$, but can have different unit form the unit of $B$. Symbols "$\twoheadrightarrow$" and "$\hookrightarrow$" denote unit preserving surjective and injective homomorphisms of algebras.
\end{proof}

Choose an orthonormal basis with respect to the form $g$ in $W_\Real$. This gives an orthonormal basis in $W$ and, in this basis, the equation of $W_1$ coincides with the equation of $S_\Complex$ given in the introduction, so $W_1 \cong S_\Complex$ and the theorem 1 follows from the following general observation:
\begin{lem}
If a finite group $\Gamma$ acts on an affine variety $X$ over $\Complex$ and $p$ is an idempotent in the group algebra of $\Gamma$ then the algebra $p \Complex[X] \ast \Gamma p$ is isomorphic to the algebra $F_\Gamma(X, \End_\Complex(\Complex \Gamma p))$ of regular $\Gamma$-equivariant maps from $X$ to $\End_\Complex(\Complex \Gamma p)$, where $\Gamma$ acts on $\End_\Complex(\Complex \Gamma p)$ by conjugation. If, moreover, $p$ is self-adjoint and $X$ is defined over $\Real$, then the involution on $\Complex[X]$ given by $f^*(x) = \overline{f(\overline{x})}$ induces an involution on $p \Complex[X] \ast \Gamma p$, and the corresponding involution on $F_\Gamma(X, \End_\Complex(\Complex \Gamma p))$ is given by
\[
f^*(x) = (f(\overline{x}))^*, \; \text{for $x\in X$, $f\in F_\Gamma(X, \End_\Complex(\Complex \Gamma p))$.}
\]
\end{lem}
For our case one should set $X=S_\Complex$ and $p=f_i$ so that $\Complex \Gamma f_i \cong V_i$.

\section{Real structure}\label{sec5}
We denote by $A$ the *-algebra $F_\Gamma(S_\Complex, \End_\Complex(V_i))$~--- the *-algebra of regular $\Gamma$-equivariant maps from $S_\Complex$ to $\End_\Complex(V_i)$ where $S_\Complex$ is an affine variety in $\Complex^3$ given by the equation $\alpha^2+\beta^2+\gamma^2=1$ and $V_i$ is a unitary irreducible representation of $\Gamma$. The involution in $A$ is given by $f^*(x) = f(\overline{x})^*$ for $x\in S_\Complex$, $f\in A$. As it was proved in the previous section $A$ is isomorphic to $e_i \Pi^\lambda e_i$ as a *-algebra.

For any point $x\in S_\Complex$ we denote by $Stab_x\subset \Gamma$ the stabilizer of $x$.
Put
\[ 
 M_x = \{m \in \End_\Complex(V_i) | m g = g m, \, \text{for any} \, g\in Stab_x\} \; \text{--- the centralizer of $Stab_x$ in $V_i$.}
\]
$M_x$ is a *-subalgebra of $\End_\Complex(V_i)$ and it is clear that for any $x\in S_\Complex$ $f(x)\in M_x$. The opposite is true:
\begin{lem}\label{lem51}
For any tuple of points $x_1$, $x_2$, \dots, $x_n$ with pairwise disjoint orbits the map from $A$ to $M_{x_1}\times M_{x_2} \times \dots \times M_{x_n}$ given by
\[
f \longmapsto (f(x_1), f(x_2), \dots, f(x_n))
\]
is surjective.
\end{lem}
\begin{proof}
For any $y_1\in M_{x_1}$, $y_2\in M_{x_2}$, \dots, $y_n\in M_{x_n}$ take any regular map from $S_\Complex$ which at points of the form $g x_j$ accept value $g y_j g^{-1}$. Clearly such a map exists because the number of points of the form $g x_j$ is finite and if $g_1 x_j = g_2 x_k$ then $j=k$ and $g_1 y_j g_1^{-1} = g_2 y_j g_2^{-1}$. Averaging this map with respect to all elements of $\Gamma$ gives an element of $A$ which maps to $(y_1, y_2, \dots, y_n)$.
\end{proof}

There is a central *-subalgebra $A_Z$ of $A$ which consists of such $f\in A$ that $f(x)$ is scalar for all $x\in S_\Complex$. The algebra $A_Z$ is the algebra of regular functions on $S_\Complex/\Gamma$. Every irreducible *-representation $\rho$ of $A$ gives a character on $A_Z$ by the Schur's lemma, thus a point $x\in S_\Complex$ such that for $f\in A_Z$ $\rho(f) = f(x)$. The point $x$ is such that 
\[
f^*(x) = \rho(f^*) = \overline{\rho(f)} = \overline{f(x)} = f^*(\overline{x}),\; \text{for any} \; f\in A_Z,
\]
so $\overline{x} = g x$, some $g\in \Gamma$.
 
\begin{prop}\label{prop51}
If $x\in S_\Complex$ and $\overline{x} = g x$ for $g\in \Gamma$ then we are in the one of two cases:
\begin{enumerate}
\item $x = \overline{x}$, $x\in S_\Real$.
\item $x$ and $\overline{x}$ are linearly independent over $\Complex$, $Stab_x = \{-1, 1\}$ and $g$ is an element of order $4$ in $SU(2, \Complex)$.
\end{enumerate}
\end{prop}
\begin{proof}
 If we suppose that $x$ and $\overline{x}$ are linearly dependent, i.e. $x = c \overline{x}$ for $c\in \Complex$ then expressing $x$ in coordinates $x = (x_1, x_2, x_3)$ gives 
 \[
 1 = x_1^2 + x_2^2 + x_3^2 = c^2 (\overline{x_1^2} + \overline{x_2^2} + \overline{x_3^2}) = c^2.
 \]
 If $c=-1$ then $x_j^2 = -x_j \overline{x_j} \le 0$~--- contradiction, so $c = 1$ and $x = \overline{x}$. Suppose that $x$ and $\overline{x}$ are linearly independent over $\Complex$. We obtain, since $g$ is defined over $\Real$ $x = g \overline{x}$, which implies that $g^2$ is the identity on $\Complex^3$ since it has at least two eigenvalues $1$. If $Stab_x$ contains an element $h\in \Gamma$ which acts non-trivially on $S_{\Complex}$ then it's eigenspace for eigenvalue $1$ is one-dimensional, so $x$ and $\overline{x}$ cannot be linearly independent.
\end{proof}
 
 Let $m_x$ be the two-sided ideal in $A$ of maps vanishing at $x$, $m_x^* = m_x$. Clearly, $\rho$ is zero on $m_x \cap A_Z$. In general $A (m_x \cap A_Z) \neq m_x$, but the following holds:
\begin{prop}
Let $\rho$ be a *-representation of $A$ and $x\in S_\Complex$ such that $\overline{x} = g x$ . If $\rho$ vanishes on $m_x \cap A_Z$ than $\rho$ vanishes on $m_x$.
\end{prop}
\begin{proof}
Let $a\in m_x$, $a = a^*$. Consider the characteristic polynomial 
\[
p(t; X) = \det{X - a(t)} = X^k + p_{k-1}(t) X^{k-1} + \dots + p_0(t).
\]
Its coefficients $p_{j}(t)$ belong to $A_Z$ and $a(x) = 0$, so $p_j \in m_x \cap A_Z$, and it follows that $\rho(a)^k = 0$. Since $\rho(a)$ is self-adjoint $\rho(a) = 0$. Since every element $b\in m_x$ can be represented as $b = a_1 + \cplxi a_2$ whith $a_1, a_2$~--- self-adjoint elements of $m_x$ we obtain $\rho(b) = 0$.
\end{proof}
Thus every irreducible *-representation of $A$ is induced from a representation of $M_x$. Let $\rho'$ be the corresponding representation of $M_x$. Then, for any $a\in A$
\[
\rho'(a(x))^* = \rho(a)^* = \rho(a^*) = \rho'(a^*(x)) = \rho'(a(\overline{x})^*) = \rho'(a(g x)^*) = \rho'(g a(x)^* g^{-1}).
\]
In fact we have
\begin{prop} For all $a\in M_x$ $\rho'(g a g^{-1}) = \rho'(a)$. Thus $\rho'(a(x))^* = \rho'(a(x)^*)$. So every irreducible *-representation of $A$ is induced from a *-representation of $M_x$, $x\in S_\Complex$.
\end{prop}
\begin{proof}
Consider the operator $\phi$ on $M_x$ sending $a$ to $g a g^{-1}$. Since $g x = \overline{x}$ and the action of $\Gamma$ on $S_\Complex$ is defined over $\Real$ $Stab_x = Stab_{\overline{x}}$ and it follows that $g$ commutes with $Stab_x$, thus its image in $\End_\Complex(V_i)$ belongs to $M_x$. Then $\phi^2$ is the identity since $g^2 x = x$ and so $g^2 \in Stab_x$. It follows that $M_x$ can be split into two eigenspaces corresponding to eigenvalues $1$ and $-1$. Suppose there is $a\in M_x$ such that $\phi(a) = -a$. So
\[
\rho'(a) \rho'(a)^* = \rho'(a g a^* g^{-1}) = - \rho'(a a^*),
\]
but $a a^*$ has spectrum contained in $\{r\in \Real | r\ge 0\}$, so the operator on the lefthand side has spectrum contained in $\{r\in \Real | r\le 0\}$. On the other hand $ \rho'(a) \rho'(a)^* \ge 0$, so $\rho'(a) = 0$. It follows that the eigenspace with eigenvalue $-1$ belongs to the kernel of $\rho'$ which implies the statement.
\end{proof}

\begin{definition}\label{def51}
We say that $V_i$ is exceptional if there is an element of $\Gamma$ of order $4$ which acts as a scalar in $V_i$.
\end{definition}

\begin{prop}\label{prop54}
If $V_i$ is not exceptional then every irreducible *-representation of $A$ is induced from a *-representation of $M_x$ for $x\in S_\Real$.
\end{prop}
\begin{proof}
Suppose we are given an irreducible *-representation $\rho'$ of $M_x$ which induces a *-representation of $A$ and $x \neq \overline{x}$. By the proposition \ref{prop51} $\overline{x} = g x$ for $g$ of order $4$. Then, for all $a\in A$
\[
\rho'(a(x)) = \rho(a) = \rho(a^*)^* = \rho'(a^*(x))^* = \rho'(a(\overline{x})^*)^* = \rho'(a(\overline{x})) = \rho'(g a(x) g^{-1}).
\]
Since $Stab_x = \{-1, 1\}$, $M_x = \End_\Complex(V_i)$ and is simple. Thus $\rho'$ has zero kernel, which implies that $g$ commutes with elements of $\End_\Complex(V_i)$. Hence $g$ is scalar in $V_i$, so $V_i$ is exceptional.
\end{proof}

Let us introduce a $C^*$-algebra $\widetilde{A}$ of all continuous $\Gamma$-equivariant maps from $S_{\Real}$ to $\End_\Complex(V_i)$. We consider a natural map $\psi : A \To \widetilde{A}$ given by restriction from $S_\Complex$ to $S_\Real$. The map is an inclusion since $S_\Real$ is algebraically dense in $S_\Complex$. Then, any irreducible *-representation $\rho$ of $A$ by proposition \ref{prop54} is induced from a *-representation $\rho'$ of $M_x$, $x\in \Real$, which induces a *-representation $\rho_0$ of $\widetilde{A}$ such that $\rho = \rho_0 \circ \psi$. Because of the lemma \ref{lem51} we can apply a Stone-Weierstrass theorem to the image of $\psi$ to show that it is dense. Given this the theorem 2 follows.

\section{Exceptional representations} \label{sec6}
Here we are going to list all exceptional representations of $\Gamma$ according to the definition \ref{def51}.
\begin{prop}
Suppose $\Gamma$ is a finite subgroup of $SU(2, \Complex)$, $g$ be an element of order $4$ and $V$~--- an irreducible representation of $\Gamma$ such that $g$ is scalar in $V$. Then either $V$ is one-dimensional, or one of the following holds:
\begin{enumerate}
\item The $\Gamma$ is a binary dihedral group whose image in $SO(3, \Real)$ is a group of symmetries of flat polygon of $n$ number of vertices, $n$ even, $n \ge 4$. The image of $g$ is the symmetry with respect to a line passing through the center of the polygon orthogonal to the plane of the polygon. $V$ is any representation corresponding to one of the black vertices on the following picture by the McKay correspondence:
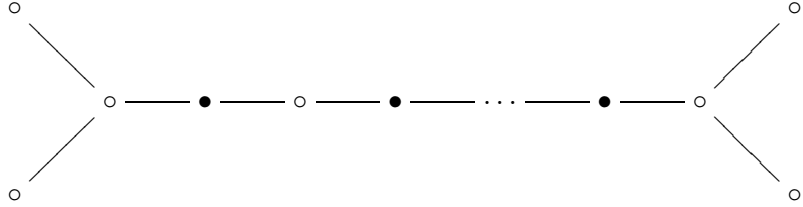
\[
\xymatrix{
\circ & & & & & & & &\circ \\
& \circ \ar@{-}[r] \ar@{-}[ul] \ar@{-}[dl] & \bullet \ar@{-}[r] & \circ \ar@{-}[r] & \bullet \ar@{-}[r] & \dots \ar@{-}[r] & \bullet \ar@{-}[r] & \circ \ar@{-}[ur] \ar@{-}[dr]& \\
\circ & & & & & & & &\circ
} 
\]
\item The $\Gamma$ is a binary octahedral group whose image in $SO(3, \Real)$ is a group of symmetries of the regular octahedron. The image of $g$ is a symmetry with respect to the line passing through opposite vertices of the octahedron. $V$ is the representation which correspond to the black vertex on the picture:
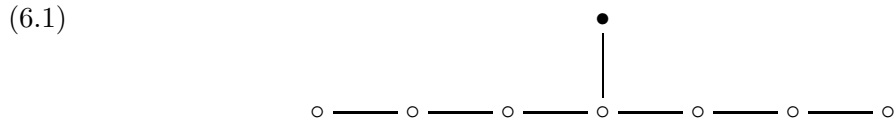
\begin{equation} \label{eq62}
\xymatrix{
& & & \bullet \ar@{-}[d] & & & \\
\circ \ar@{-}[r] & \circ \ar@{-}[r] & \circ \ar@{-}[r] & \circ \ar@{-}[r] & \circ \ar@{-}[r] & \circ \ar@{-}[r] & \circ 
}
\end{equation}
\end{enumerate}
\end{prop}
\begin{proof}
The proof consists of looking at values of characters on elements of order $4$. We consider here only the case of the binary octahedral group. Other cases are similar to this one. The McKay graph of $\Gamma$ is (\ref{eq62}). Let $g$ be an element of order $4$. Each vertex of the graph we label by the number which equals to the trace of $g$ in the corresponding representation. We collect facts we know about these numbers:
\begin{enumerate}
\item On the identity representation the trace is $1$.
\item On the tautological two-dimensional representation $V$ the trace is $0$, so for any other representation $V_i$ the trace on $V\otimes V_i$ is zero. Thus the sum of labels of neighbours of $i$ is zero.
\item For any even vertex we have $g^2$ acts identically in the corresponding representation, so its spectrum contains only $1$ and $-1$.
\item For any odd vertex we have $g^2$ acts as negative identity in the corresponding representation, so its spectrum contains only $\cplxi$ and $-\cplxi$.
\end{enumerate}
We obtain only two possible labelings which satisfy the conditions above:
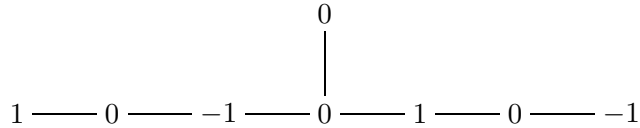
\[
\xymatrix{
& & & 0 \ar@{-}[d] & & & \\
1 \ar@{-}[r] & 0 \ar@{-}[r] & -1 \ar@{-}[r] & 0 \ar@{-}[r] & 1 \ar@{-}[r] & 0 \ar@{-}[r] & -1}
\]
and
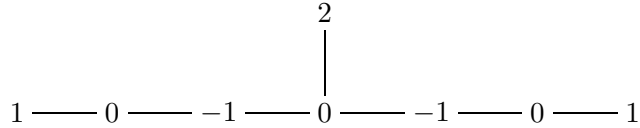
\[
\xymatrix{
& & & 2 \ar@{-}[d] & & & \\
1 \ar@{-}[r] & 0 \ar@{-}[r] & -1 \ar@{-}[r] & 0 \ar@{-}[r] & -1 \ar@{-}[r] & 0 \ar@{-}[r] & 1} 
\]
The first case cannot give an exceptional case of dimension greater then $1$ since absolute value of all traces is not more than one. The second case gives an exceptional case of dimension $2$ if there exists an element of $\Gamma$ with such traces. The sum of squares of absolute value is $4$ in the first case and $8$ in the second. This number equals to the cardinality of the centralizer of $g$. Then, considering a homomorphism $SU(2, \Complex) \To SO(3, \Real)$, if $g$ is a preimage of the symmetry with respect to a line passing through two opposite vertices of the octahedron, it commutes with preimages of all rotations fixing this line, which number is 4. We see that centralizer of $g$ has cardinality $8$, so its traces are given by the second labeling. The first labeling corresponds to the symmetries with respect to lines, passing through the middlepoints of edges of the octahedron.
\end{proof}

\section{Trivializing bundles} \label{sec7}
In this section we are going to prove the theorem 3.
Suppose we have a finite group $\Gamma \subset SU(2, \Complex)$ containing $-I$ and its irreducible unitary representation $V_i$. Then the group $\Gamma' = \Gamma/\{-I, I\}$ is a subgroup of $SO(3, \Real)$ and acts on $\End_\Complex(V_i)$ by conjugation. Let $\cA$ be the *-algebra of continuous maps $f$ from the unit sphere $S_\Real$ to $\End_\Complex(V_i)$ such that $f(g x) = g f(x) g^{-1}$ for all $g\in \Gamma'$. Note that $S_\Real/\Gamma'$ in all cases is homeomorphic to $S_\Real$. If $\Gamma$ is cyclic then $V_i$ is one-dimensional, so 
\begin{prop}
If $\Gamma$ is cyclic then $\cA$ is isomorphic to the *-algebra
 of continuous functions $S_\Real \To \Complex$. 
\end{prop}
Suppose that $\Gamma$ is not cyclic. Then it is a binary dihedral, tetrahedral, octahedral or icosahedral group so that $\Gamma'$ is the usual dihedral, tetrahedral or icosahedral group. In all cases there is a fundamental domain for $\Gamma'$ on $S_\Real$ of the form:
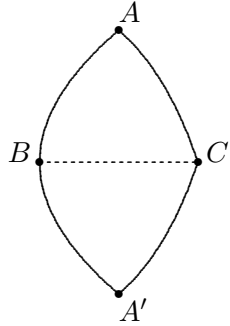
\[
\begin{picture}(100,115)(0,-5)
\curvedashes{0,1,2}
\tagcurve(10,50,20,50,80,50,90,50)
\curvedashes{}
\tagcurve(90,110, 50,100, 20,50, 50,0, 90,-10)
\tagcurve(10,110, 50,100, 80,50, 80,10)
\tagcurve(10,-10, 50,0, 80,50, 80,90)
\put(20,50){\circle*{3}}
\put(50,100){\circle*{3}}
\put(50,0){\circle*{3}}
\put(80,50){\circle*{3}}
\put(8,50){$B$}
\put(50,103){$A$}
\put(83,50){$C$}
\put(50,-10){$A'$}
\end{picture}
\]
where $A$ and $A'$ are centers of faces, $B$ is the center of an edge and $C$ is a vertex of the correspoding dihedra, tetrahedra, octahedra or icosahedra. Denote this fundamental domain by $\cF$. The group is generated by three elements $a$, $b$, $c$, where $a$ is a rotation around $A$, $b$ is a rotation around $B$ and $c$ is a rotation around $C$. Elements $a$, $b$, $c$ can be choosen in such a way that $c A = A'$, $b A' = A$ and $abc = e$, $e$ denotes the identity. Then the stabilizer of $A$ is generated by $a$, the stabilizer of $B$ by $b$, the stabilizer of $C$ by $c$ and the stabilizer of $A'$ by $c a c^{-1}$. Other points of the fundamental domain have trivial stabilizers. Every orbit of the group intersects the fundamental domain in exactly one point except orbits passing through the boundary, i.e. the element $b$ maps points of the segment $B A'$ to $B A$ and $c$ maps points of $C A$ to $C A'$.
It follows that the quotient $S_\Real/\Gamma'$ can be obtained by gluing the fundamental domain along actions of $b$ and $c$. 
We are going to consider the following class of $C^*$-algebras.
\begin{definition} \label{def71}
Suppose we are given the following data:
\begin{enumerate}
\item a finite dimensional hermitian vector space $H$,
\item a continuous map $m_b: BA' \To SU(H)$ which is constant in the neighbourhoods of the endpoints of $BA'$,
\item a continuous map $m_c: C A \To SU(H)$ which is constant in the neighbourhoods of the endpoints of $CA$,
\item a *-subalgebra $M_P \subset \End(H)$ for every point $P$ in the set $\{A, B, C, A'\}$,
\end{enumerate}
which satisfy the following conditions:
\begin{enumerate}
\item $m_b(B)$ commutes with $M_B$,
\item $m_c(C)$ commutes with $M_C$,
\item $m_b(A') m_c(A)$ commutes with $M_A$,
\item $m_c(A) M_A m_c(A)^{-1} = M_{A'}$.
\end{enumerate}
Then we denote by $\cC(H, m_b, m_c, (M_P))$ the $C^*$-algebra of continuous maps $f:\cF \To \End(H)$ such that:
\begin{enumerate}
\item $f(P)\in M_P$ for $P\in \{A, B, C, A'\}$,
\item $f(bx) = m_b(x) f(x) m_b(x)^{-1}$ for $x\in BA'$,
\item $f(cx) = m_c(x) f(x) m_c(x)^{-1}$ for $x\in CA$.
\end{enumerate}
\end{definition}

It is clear that the algebra $\cA$ is isomorphic to the algebra $\cC(V_i, \widetilde{b}, \widetilde{c}, (M_P))$, where $\widetilde{b}$ and $\widetilde{c}$ are constant maps with \[
\widetilde{b} = \frac{V_i(b)}{\sqrt[d]{\det(V_i(b))}}, \;
\widetilde{c} = \frac{V_i(c)}{\sqrt[d]{\det(V_i(c))}},
\]
here $V_i(b)$ denotes the action of $b$ on $V_i$ and the same for $V_i(c)$.

Suppose we are given two algebras $\cC^1 = \cC(H, m_b^1, m_c^1, (M^1_P))$ and $\cC^2 = \cC(H, m_b^2, m_c^2, (M^2_P))$ with the same space $H$.
\begin{definition} \label{def72}
A map $t: \cF \To SU(H)$ is called a morphism between $\cC^1$ and $\cC^2$ if it satisfies the following properties:
\begin{enumerate}
\item $t$ is continuous in all points except, possibly, $A$, $B$, $C$ and $A'$.
\item If $P$ is any point among $A$, $B$, $C$ and $A'$ then there is a neighbourhood of $P$ such that for any point $x$ from the neighbourhood and $u\in M^1_P$ 
\begin{equation}
t(x) u t(x)^{-1} = t(P) u t(P)^{-1} \in M^2_P. 
\end{equation}
\item $m_b^2(x) t(x) = t(bx) m_b^1(x)$ for all $x\in BA'$, except, possibly, $B$ and $A'$.
\item $m_c^2(x) t(x) = t(cx) m_c^1(x)$ for all $x\in CA$, except, possibly, $C$ and $A$.
\end{enumerate}
\end{definition}

\begin{prop} \label{prop72}
If $t$ is a morphism then there is a homomorphism of $C^*$-algebras $\widetilde{t}: \cC^1 \To \cC^2$ defined by
\[
(\widetilde{t} f) (x) = t(x) f(x) t(x)^{-1},\; \text{for all $x\in \cF$, $f\in \cC^1$.}
\]
\end{prop}
\begin{proof}
First we prove that continuous maps are mapped to continuous maps. This property is obvious in all points except $A, B, C, A'$. If $P$ is one of $A$, $B$, $C$, $A'$ and $f$ is an element of $\cC$, then 
\[
\lim_{x\To P} t(x) f(x) t(x)^{-1} - t(P) f(P) t(P)^{-1} = \lim_{x\To P} t(x) (f(x) - f(P)) t(x)^{-1} = 0,
\]
the first equality follows from the condition (2) of the definition \ref{def72} and the last equality is true because $t(x)$ and $t(x)^{-1}$ are bounded.
To prove the statement we must choose any $f\in \cC^1$ and show that $g = \widetilde{t} f \in \cC^2$. We check conditions of the definition \ref{def71} one by one.
\begin{enumerate}
\item For $P\in \{A, B, C, A'\}$ $g(P) = t(P) f(P) t(P)^{-1} \in M_P^{2}$ since $f(P) \in M_P^1$ by the condition (2) of the definition \ref{def72}.
\item For $x \in BA'$, $x \neq B$, $x \neq A'$ 
\[
\begin{split}
g(b x) m_b^2(x) = t(b x) f(b x) t(b x)^{-1} m_b^2(x) = t(b x) f(b x) m_b^1(x) t(x)^{-1} 
\\= t(b x) m_b^1(x) f(x) t(x)^{-1} = m_b^2(x) t(x) f(x) t(x)^{-1} = m_b^2(x) g(x).
\end{split}
\]
Since $g$ is already known to be continuous the statement follows for $B$ and $A'$.
\item Can be proved analogously to (2).
\end{enumerate}
\end{proof}

Using morphisms we can show that
\begin{thm}
If $\cC^1 = \cC(H, m_b^1, m_c^1, (M^1_P))$ and $\cC^2 = \cC(H, m_b^2, m_c^2, (M^2_P))$ are such that for $P = A, B, C$ the *-algebra $M^1_P$ is isomorphic to $M^2_P$ then $\cC^1$ and $\cC^2$ are isomorphic $C^*$-algebras. 
\end{thm}
\begin{proof}
Choose $u_P\in SU(H)$ such that $M^2_P = u_P M^1_P u_P^{-1}$ for $P\in \{A, B, C\}$.  We first define $t(P) = u_P$, $t(A') = m_c^2(A) u_A m_c^1(A)^{-1}$. Then we choose $t$ on $AB$ such that $t$ equals to $u_A$ in a neighbourhood of $A$, to $u_B$ in a neighbourhood of $B$ and is continuous. Since $SU(H)$ is linearly connected it is possible to do so. Then for $AC$ we do analogously. For $BA'$ and $CA'$ (except $B$, $C$, $A'$) we put
\[
t(x) = m_b^2(x)^{-1} t(b x) m_b^1(x) \; \text{and} \; t(x) = m_c^2(c^{-1} x) t(c^{-1} x) m_c^1(c^{-1} x)^{-1} \; \text{correspondingly.}
\]
In such a way we achieve properties (3) and (4) of the definition \ref{def72}. 

Next we consider small disks around points $A$, $B$, $C$, $A'$. For small disk around $A$ we put $t(x) = u_A$ and the condition (2) will be satisfied. For small disk around $B$ we already defined $t$ on points of segment $BA$ to be $u_B$ and on points of segment $BA'$ to be $m_b^2(B)^{-1} u_B m_b^1(B)$. Since $m_b^1(B)$ commutes with elements of $M_B^1$ and $m_b^2(B)^{-1}$ commutes with elements of $M_B^2 = u_B M_B^1 u_B^{-1}$ we get $u_B^{-1} m_b^2(B)^{-1} u_B$ and, hence $u_B^{-1} m_b^2(B)^{-1} u_B m_b^1(B)$ commute with elements of $M_B^1$. Then there exists a path $p: [0, 1] \To SU(H)$ connecting $I_H$ and $u_B^{-1} m_b^2(B)^{-1} u_B m_b^1(B)$ such that every point of the path commutes with elements of $M_B^1$. Take the path $u_B p$. It connects $u_B$ with $m_b^2(B)^{-1} u_B m_b^1(B)$, so we can use this path to define $t$ on arcs of circles centered at $B$ with endpoints on $BA$ and $BA'$. Thus we obtain a map defined in a disk centered at $B$, which is continuous at all points except $B$ and for $u\in M_B^1$, $x$ in this disk
\[
t(x) u t(x)^{-1} = u_B p(\tau) u p(\tau) u_B^{-1} = u_B u u_B^{-1} = t(B) u t(B)^{-1}, \; \text{for some $\tau \in [0, 1]$}.
\]
The same can be done for $C$ and $A'$ since for $C$:
\[
u_C^{-1} m_c^2(C) t(C) m_c^1(C)^{-1} = u_C^{-1} m_c^2(C) u_C m_c^1(C)^{-1}
\]
commutes with elements of $M_C^1$. In the neighbourhood of $A'$ we have value $m_b^2(A')^{-1} u_A m_b^1(A')$ for points of $BA'$ and $m_c^2(A) u_A m_c^1(A)^{-1}$ for points of $CA'$ and for $A'$ itself. So we check that
\[
\begin{split}
(m_c^2(A) u_A m_c^1(A)^{-1})^{-1} m_b^2(A')^{-1} u_A m_b^1(A') = 
m_c^1(A) u_A^{-1} m_c^2(A)^{-1} m_b^2(A')^{-1} u_A m_b^1(A') \\= 
m_c^1(A) \left[ u_A^{-1} m_c^2(A)^{-1} m_b^2(A')^{-1} u_A m_b^1(A') m_c^1(A) \right] m_c^1(A)^{-1}
\end{split}
\]
commutes with elements of $M^1_{A'}$ since $M^1_{A'} = m_c^1(A) M^1_A m_c^1(A)^{-1}$ and the element in square brackets commutes with elements of $M^1_A$. Indeed, $m_b^1(A') m_c^1(A)$ commutes with elements of $M^1_A$ because of the property (3) of the definition \ref{def71}, and $m_c^2(A)^{-1} m_b^2(A')^{-1}$ commutes with elements of $M^2_A$ by the same reason. Using $M^1_A = u_A^{-1} M^2_A u_A$ we obtain the required property.

Then using the fact that $SU(H)$ is simply connected we can extend $t$ to the whole of $\cF$ which will give a morphism in the sense of the definition \ref{def72}. By the proposition \ref{prop72} we obtain a homomorphism of $C^*$-algebras. Taking the inverse of $t$ at each point of $\cF$ gives an inverse homomorphism. Thus the theorem is proved.
\end{proof}

The theorem 3 is a trivial consequence of the theorem above since both algebras in the statement of the theorem 3 can be represented according to the definition \ref{def71}, they have the same space $H = V_i = \Complex^d$ and isomorphic subalgebras $M_P$ at points $A$, $B$, $C$.

\end{document}